\documentclass[]{interact}

\usepackage{epstopdf}
\usepackage[caption=false]{subfig}

\usepackage{tikz}
\usepackage{amsmath}
\usepackage{hyperref}

\usepackage[numbers,sort&compress]{natbib}
\bibpunct[, ]{[}{]}{,}{n}{,}{,}
\makeatletter
\def\NAT@def@citea{\def\@citea{\NAT@separator}}
\makeatother

\numberwithin{equation}{section}
\newenvironment{rcases}
  {\left.\begin{aligned}}
  {\end{aligned}\right\rbrace}
\newcommand{\eps}{\varepsilon}
\newcommand*{\R}{\mathbb{R}}

\renewcommand*{\S}{\mathcal{S}}
\newcommand*{\cc}{\mathbb{C}}
\renewcommand*{\phi}{\varphi}
\DeclareMathOperator{\Imag}{Im}
\newcommand*{\X}{\mathcal{X}}

\theoremstyle{plain}

\theoremstyle{definition}

\theoremstyle{remark}

\begin{document}


\articletype{}

\title{On the open sea propagation of two-dimensional rotational water waves generated by a moving bed}

\author{
\name{Frederick Moscatelli\textsuperscript{a}\thanks{CONTACT Frederick Moscatelli. Email: frederick.moscatelli@univie.ac.at}}
\affil{\textsuperscript{a}Faculty of Mathematics, University of Vienna\\
Oskar-Morgenstern-Platz 1, 1090 Vienna, Austria}
}

\maketitle

\begin{abstract}
We study the propagation of two-dimensional tsunami waves triggered by a seaquake in the open sea in the presence of underlying wind-generated currents, 
corresponding to background flows of constant vorticity. A suitable scaling of the governing equations introduces dimensionless parameters, of particular interest being 
the setting of linear waves that only depend on the vertical movement of the sea bed. We use Fourier analysis methods to extract formulae for the function $f$ 
which describes the vertical displacement of the water's free surface. We show that the results are particularly useful in the physically relevant shallow-water 
regime: in the irrotational case the predictions fit well with the observed behaviour of some historical tsunamis. In other situations, the stationary-phase principle gives insight into the asymptotic behaviour of $f$.
\end{abstract}


\begin{amscode}
45G99, 58J32, 76B03
\end{amscode}

\section{Introduction}

The aim of this paper is to pursue a mathematical analysis of tsunami waves triggered by an earthquake in the open sea. 

Tsunamis are water waves which are generated by an impulsive disturbance that vertically displaces the water column, 
most commonly by tectonic activity, but sometimes by other means as well, such as underwater landslides or asteroid impact. 
Hundreds of tsunamis were confirmed since 1900, about 80\% of which occurred in the Pacific Ocean. Most of these tsunami waves 
have small amplitude, being nondestructive and barely detectable in the ocean (less than 0.3 m high, the average wave height in 
the Pacific Ocean being about 2.5 m). Large, destructive tsunamis, causing loss of life and major coastal
destruction, take place typically once every several decades. The highest tsunami wave that was ever measured exceeded 500 m and 
occurred on 9 July 1958, in the Lituya Bay inlet of the Gulf of Alaska, caused by a massive rock dislocated by an earthquake and 
plunging 900 m down into the water (see \cite{Con11}). Let us also mention that 66 million years ago, an asteroid with a diameter of about 10 km struck the 
waters of the Gulf of Mexico near the Yucat\'an Peninsula, and it is estimated that the impact generated tsunami 
waves more than 1.5 km high that caused sea levels to rise in all corners of the Earth (see \cite{r}); the subsequent cataclysmic events -- 
plumes of aerosol, soot and dust filling the air, and wildfires starting as flaming pieces of material blasted from the impact re-entered the atmosphere and rained down -- 
ended the era of the dinosaurs and triggered a mass extinction of about 75\% of animal and plant life on Earth. Since it appears \cite{as} 
that up to 75\% of all recorded tsunamis were generated by undersea earthquakes, in this paper we focus on this type. Three events stand out in the last 
100 years: the 22 May 1960 Chile tsunami (caused by the largest earthquake ever recorded and with waves propagating across the Pacific Ocean, causing havoc in 
Hawaii and in Japan), the 26 December 2004 tsunami (that killed more than 200.000 people around the shores of the Indian Ocean), and the tsunami off the
coast of Japan on 11 March 2011 (that killed thousands of people and triggered a nuclear accident). In all three cases the earthquake was undersea, along a fault line, in 
which case the generated waves are typically two-dimensional and with long wavelengths, of the order of 100 km (see the discussions in \cite{as,Con11,dd}). Given that the 
average ocean depth is about 4 km, this means that tsunami waves in the open ocean are shallow water waves. While near the shore, where the water depth gradually diminishes, 
nonlinear effects become dominant (and induce wave-breaking), a long-standing issue was whether in the open sea the dominant linear behaviour of the flow is greatly affected by the 
accumulation of nonlinear factors. Some authors (see \cite{cr, la, se}) advocated that KdV-like models apply, but even in the case of propagation over the maximal distance that 
is possible (across the Pacific Ocean, which was the case for the 1960 Chile tsunami), the scale analysis reveals that weakly nonlinear theory is not adequate (see the 
discussion in \cite{as, Con09, Con11, cj, st}): the effects remain negligible in the open sea. Advances in the modelling of the 
tsunami wave propagation in the open sea, triggered by an undersea earthquake, were obtained in \cite{CG12, dd} for the setting of a still sea prior to the perturbation. 
In this paper we extend these considerations to accommodate the presence of underlying currents. These are typically modelled by constant-vorticity flows 
(see the discussion in \cite{TDSP88, Ewi90}) but even the irrotational setting (flows with zero vorticity) might admit uniform tidal currents.

\section{The governing equations}

\begin{figure}[!ht]
\centering
\begin{tikzpicture}
\draw[blue] (-5,3) -- (5,3);
\draw[dashed] (0,0) -- (0,3);
\draw[<->] (-3,0) -- (-3,3) node[anchor=north east] {$d$};
\foreach \x in {0.5,1,1.5,2,2.5}
    \draw[thick,->] (0,\x) -- (0.5 + \x/2,\x);
\draw[thick] (-5,0) -- (5,0);
\draw[color = red, thick,->] (0,3) -- (2,3) node[anchor=south] {$A d + B$};
\draw[color=red,thick,->] (0,0) -- (0.5,0) node[anchor=north east] {$B$};
\foreach \x in {-5,-4,-3,-2,-1,0,1,2,3,4}
    \draw[thick] (\x,0) -- (\x + 0.1,-0.2);
\foreach \x in {-5,-4,-3,-2,-1,0,1,2,3,4}
    \draw[thick] (\x+0.5,0) -- (\x + 0.6,-0.2);
\end{tikzpicture}
\label{Model at time 0}
\caption{Pure current flow (in the absence of waves) at the initial time $T=0$}
\end{figure}
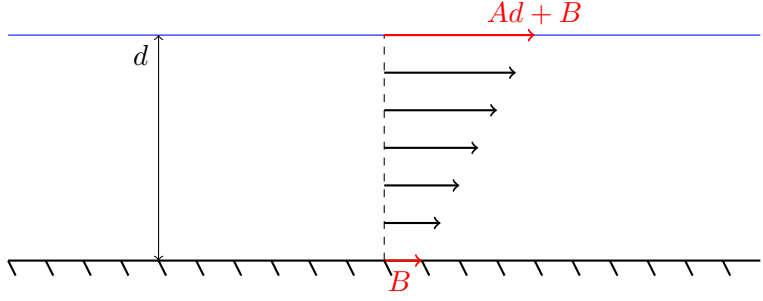

We consider an analogous setup as in \cite{CG12}: we consider two-dimensional surface waves with the unbounded $X$-direction corresponding to 
the direction of wave propagation. The water's free surface is $Y = d + F(X,T)$, where $d$ is the average depth of the sea, the seabed being 
$Y = H(X,T)$. At time $T=0$, we assume that $F(X,0)=0$ and $H(X,0)=0$. Then the seabed moves for a short time in some compact region 
$X \in [-L,L]$ and is flat after that again. We are ultimately interested in the deformation of the free surface, i.e. in $F(X,T)$. To do this, we study a 
partial differential equation that describes its evolution in time. Consider $(U,V)$ the velocity field of the flow in rectangular Cartesian coordinates $(X,Y)$. We assume that the water is inviscid and that the resulting flow is with constant vorticity. Note that the presence of non-uniform currents makes the hypothesis of irrotational flow inappropriate and, since in a setting in which the waves are long compared to the water depth, the existence of a non-zero mean vorticity is important rather than its specific distribution \cite[cf.][]{TDSP88}. This is not merely a mathematical simplification: wind-generated currents in shallow waters with nearly flat beds have been shown to be accurately described as flows with constant vorticity \cite{Ewi90}. We now present the governing equations for water flows \cite{Con11}. We have the incompressible Euler equations
\begin{align}\label{Euler PDE}
\begin{rcases}
    U_X + V_Y &= 0,\\
    \rho (U_T + U U_X + V U_Y) &= - P_X,\\
    \rho (V_T + U V_X + V V_Y) &= - P_Y - \rho g,
\end{rcases}
\end{align}
in $H(X,T) < Y < d +  F(X,T)$, where $\rho$ is the constant density, $g$ is the constant acceleration of gravity and $P$ is the pressure. The effects of surface tension are negligible for wavelengths greater than a few centimetres \cite{BK75,Lig96}, hence the major factor governing the wave motion is the balance between gravity and the inertia of the system. We 
also have the kinematic boundary conditions
\begin{align}\label{bc on surface, physical}
    V = F_T + U F_X \quad \text{on the free surface }Y=d + F(X,T)
\end{align}
and
\begin{align}\label{bc on seabed, physical}
    V = H_T + U H_X \quad \text{on the rigid bed }Y=H(X,T),
\end{align}
which ensure that particles on these boundaries are confined to them at all times. We also have the dynamic boundary condition
\begin{align}\label{bc pressure, physical}
    P - P_\text{atm} = 0\quad \text{on the free surface }Y=d + F(X,T),
\end{align}
where $P_{\text{atm}}$ is the constant atmospheric pressure. This decouples the motion of the water from the motion of the air above it \cite{Joh97}. The system \eqref{Euler PDE}-\eqref{bc pressure, physical} is to be solved with the following initial conditions
\begin{align}\label{ic, physical}
    F(X,0)=0,\quad U(X,Y,0) = A Y + B,\quad V(X,Y,0) = 0.
\end{align}
These initial conditions express the assumption that initially, at time $T=0$, there is only a horizontal flow whose speed depends linearly on the depth. This is a generalization of \cite{CG12}, where it was assumed that the water is completely at rest at $T=0$, a case which corresponds to $A = B = 0$. A two-dimensional water flow that has constant vorticity at some instant will have this feature at all other times (see the discussion in Chapter 1 of \cite{Con11}). By the initial condition \eqref{ic, physical} this implies
\begin{align}\label{vorticity, physical}
    U_Y - V_X = A
\end{align}
for $H(X,T) < Y < d +  F(X,T)$. Fig \ref{Model at time 0} illustrates our set-up in the case $A,B >0$.

The first step in the analysis of the system \eqref{Euler PDE}-\eqref{vorticity, physical} will be to non-dimensionalize and scale the equations. 
First, it will be convenient to introduce the non-dimensional excess pressure $p$ relative to the hydrostatic pressure distribution by
\begin{align}\label{def of p}
    P = P_\text{atm} + \rho g (d-Y) + \rho g d p.
\end{align}
To obtain meaningful scales, we introduce the average or typical wavelength $\lambda$ of the wave, we use $\sqrt{ g d }$ as a scale for the wave speed and $\lambda/\sqrt{g d}$ as a time scale. We use the change of variables
\begin{align}\label{first scaling 1}
    X = \lambda x,\quad Y = d y,\quad T = \frac{\lambda}{\sqrt{gd}}t,\quad U = \sqrt{g d}u,\quad V = \frac{d\sqrt{gd}}{\lambda}v
\end{align}
with 
\begin{align}\label{first scaling 2}
    F = a f,\quad H = a h,
\end{align}
where $a$ is a typical, perhaps maximal, amplitude of the wave. Note that \eqref{first scaling 2} should be interpreted as ensuring 
that the variations of the wave and of the seabed are of comparable size.

We introduce the parameters
\begin{align}\label{def of eps and delta}
    \eps = \frac{a}{d},\quad \delta = \frac{d}{\lambda}\,,
\end{align}
where $\eps$ measures the relative size of the amplitude to the average water depth and $\delta$ measures the average water depth to the wavelength. 
Now we get the non-dimensional equations
\begin{align}\label{non-dim orig}
\begin{rcases}
&\begin{cases}
    u_x + v_y = 0\\
    u_t + u u_x + v u_y = -p_x\\
    \delta^2(v_t + u v_x + v v_y) = -p_y\\
    u_y - \delta^2 v_x = \sqrt{\frac{d}{g}}A
\end{cases}\quad \text{in }\varepsilon h(x,t) < y < 1 + \varepsilon f(x,t)\\
    &p = \varepsilon f \qquad \qquad \quad \text{on }y = 1 + \varepsilon f(x,t)\\
    &v =\varepsilon(f_t + u f_x)\quad \text{on }y = 1 + \varepsilon f(x,t)\\
    &v =\varepsilon(h_t + u h_x) \quad \text{on }y=\varepsilon h(x,t)\\
    &f(x,0)=0,\quad u(x,y,0) =\sqrt{\frac{d}{g}}A y + \frac{B}{\sqrt{gd}},\quad v(x,y,0) = 0
\end{rcases}.
\end{align}
The magnitudes of $\eps$ and $\delta$ correspond to the different general types of water wave problem: The limits $\delta \to 0$ and $\delta \to \infty$ produce the 
shallow water and deep water regime, respectively. On the other hand, $\eps \to 0$ corresponds to regime of waves of small amplitude (see the discussion in \cite{Joh97,CJ08}).

We want to linearize the problem and in this case we will want to let $\eps \to 0$ and keep $\delta$ fixed. But the system \eqref{non-dim orig} shows that $v$ and $p$ are of order $\eps$ and hence also $u$. Hence we use the following scaling:
\begin{align}\label{second scaling}
    u \mapsto \sqrt{\frac{d}{g}}Ay + \frac{B}{\sqrt{g d}} +\varepsilon u, \quad v \mapsto \varepsilon v, \quad p\mapsto \varepsilon p,
\end{align}
where we avoid a new notation. We then obtain the system 
\begin{align}\label{non-dim new 1}
    \begin{rcases}
    u_x + v_y &= 0\\
    u_t + \left(\sqrt{\frac{d}{g}}A y + \frac{B}{\sqrt{gd}} + \varepsilon u\right)u_x + v\left(\sqrt{\frac{d}{g}}A + \varepsilon u_y\right) &= - p_x\\
    \delta^2 \left(v_t + \left(\sqrt{\frac{d}{g}}A y + \frac{B}{\sqrt{gd}} + \varepsilon u\right)v_x + \varepsilon v v_y\right) &= -p_y\\
    u_y - \delta^2 v_x &= 0
    \end{rcases}
\end{align}
in $\varepsilon h(x,t) < y < 1 + \varepsilon f(x,t)$ and 
\begin{align}\label{non-dim new 2}
    \begin{rcases}
    &p = f \quad \text{on }y = 1 + \varepsilon f(x,t)\\
    &v = f_t + \left(\sqrt{\frac{d}{g}}A (1 + \varepsilon f(x,t)) + \frac{B}{\sqrt{gd}} + \varepsilon u\right)f_x\quad \text{on }y = 1 + \varepsilon f(x,t)\\
    &v= h_t + \left(\sqrt{\frac{d}{g}}A \varepsilon h(x,t) + \frac{B}{\sqrt{gd}} + \varepsilon u\right)h_x \quad \text{on }y=\varepsilon h(x,t)\\
    &f(x,0) = 0,\quad u(x,y,0) = 0, \quad v(x,y,0) = 0
    \end{rcases}.
\end{align}
We now linearize this system by taking $\eps \to 0$ and obtain
\begin{align}\label{reduced eq}
    \begin{rcases}
    &\begin{cases}
    u_x + v_y = 0\\
    u_t + \left(\sqrt{\frac{d}{g}}A y + \frac{B}{\sqrt{gd}}\right)u_x + \sqrt{\frac{d}{g}}A v = -p_x\\
    \delta^2\left(v_t + \left(\sqrt{\frac{d}{g}}A y + \frac{B}{\sqrt{gd}}\right)v_x\right) = -p_y\\
    u_y - \delta^2 v_x = 0
    \end{cases}\quad \text{in }0 < y < 1\\
    &p = f \quad \text{on }y = 1\\
    &v = f_t + \left(\sqrt{\frac{d}{g}}A + \frac{B}{\sqrt{gd}}\right)f_x \quad \text{on }y=1\\
    &v = h_t + \frac{B}{\sqrt{gd}} h_x \quad \text{on }y=0\\
    &f(x,0)=0,\quad u(x,y,0) = 0,\quad v(x,y,0)=0
    \end{rcases}.
\end{align}
By the first equation in \eqref{reduced eq}, there exists a stream function $\psi$ such that 
\begin{align}\label{def psi}
    u = \psi_y,\quad v = -\psi_x \quad \text{in }0<y<1.
\end{align}
If we plug this in \eqref{reduced eq}, we get 
\begin{align}\label{reduced eq psi}
    \begin{rcases}
    &\begin{cases}
    \psi_{yy} + \delta^2 \psi_{xx} = 0\\
    \psi_{yt} + \left(\sqrt{\frac{d}{g}}A y + \frac{B}{\sqrt{gd}}\right)\psi_{xy} - \sqrt{\frac{d}{g}}A \psi_{x} = -p_x\\
    \delta^2\left(\psi_{x t} + \left(\sqrt{\frac{d}{g}}A y + \frac{B}{\sqrt{gd}}\right)\psi_{x x}\right) = p_y
    \end{cases}\quad \text{in }0 < y <1 \\
    &p = f \quad \text{on }y = 1\\
    &\psi_x = -f_t - \left(\sqrt{\frac{d}{g}}A + \frac{B}{\sqrt{gd}}\right)f_x \quad \text{on }y=1\\
    &\psi_x = -h_t - \frac{B}{\sqrt{gd}} h_x \quad \text{on }y=0\\
    &f(x,0)=0,\quad \psi_y(x,y,0) = 0,\quad \psi_x(x,y,0)=0
    \end{rcases}.
\end{align}
Note that we want to eventually find $f$. But $p$ and $\psi$ are also unknown, only $h$ is given. We will now derive equations only involving $\psi$ and $h$.

We first differentiate the second equation in \eqref{reduced eq psi} with respect to $x$ and $t$, respectively, to obtain
\begin{align}\label{second x}
    \psi_{xyt} + \left(\sqrt{\frac{d}{g}}A y + \frac{B}{\sqrt{gd}}\right)\psi_{xxy} - \sqrt{\frac{d}{g}}A \psi_{xx} = -p_{xx}
\end{align}
and 
\begin{align}\label{second t}
    \psi_{ytt} + \left(\sqrt{\frac{d}{g}}A y + \frac{B}{\sqrt{gd}}\right)\psi_{xyt} - \sqrt{\frac{d}{g}}A \psi_{xt} = -p_{xt}.
\end{align}
Furthermore, we can differentiate the fifth equation in \eqref{reduced eq psi} with respect to $x$ to obtain
\begin{align}\label{fifth x}
    \psi_{xx} &= -f_{x t} - \left(\sqrt{\frac{d}{g}}A + \frac{B}{\sqrt{gd}}\right)f_{xx} \quad \text{on }y=1.
\end{align}
Observe that the fourth equation in \eqref{reduced eq psi} implies that we can exchange the derivatives of $f$ in \eqref{fifth x} by the corresponding derivatives of $p$. But then we can insert the equations \eqref{second x} and \eqref{second t} in \eqref{fifth x} to obtain
\begin{align*}
    \psi_{xx} &=\psi_{y t t} + 2 \left(\sqrt{\frac{d}{g}}A + \frac{B}{\sqrt{gd}}\right)\psi_{xyt}+ \left(\sqrt{\frac{d}{g}}A + \frac{B}{\sqrt{gd}}\right)^2\psi_{x x y} \\&- \sqrt{\frac{d}{g}}A\left(\sqrt{\frac{d}{g}}A + \frac{B}{\sqrt{gd}}\right)\psi_{x x} - \sqrt{\frac{d}{g}}A \psi_{x t}
\end{align*}
on $y = 1$. We introduce the constant
\begin{align}\label{def C}
    C = \sqrt{\frac{d}{g}}A + \frac{B}{\sqrt{gd}},
\end{align}
which is the non-dimensional horizontal velocity of the flow at the surface. This allows us to rewrite the previous equation as 
\begin{align*}
    \left(1 + \sqrt{\frac{d}{g}}A C\right)\psi_{xx} = \psi_{y t t} + 2 C \psi_{x y t} + C^2 \psi_{x x y} - \sqrt{\frac{d}{g}}A \psi_{x t}
\end{align*}
on $y = 1$. By observing the quadratic expression on the right-hand side and introducing the differential operator
\begin{align}\label{def T}
    S = C \partial_x + \partial_t,
\end{align}
we can rewrite this again as 
\begin{align*}
    \left(1 + \sqrt{\frac{d}{g}}A C\right)\psi_{xx} &= S^2 \psi_y - \sqrt{\frac{d}{g}}A \psi_{x t}.
\end{align*}
So we obtain
\begin{align}\label{reduced eq final}
    \begin{rcases}
    \psi_{yy} + \delta^2\psi_{xx} &= 0,\quad \text{in } 0< y < 1\\
    \left(1 + \sqrt{\frac{d}{g}}A C\right)\psi_{xx} &= S^2 \psi_y - \sqrt{\frac{d}{g}}A \psi_{x t},\quad \text{on }y=1\\
    \psi_x &= -h_t - \frac{B}{\sqrt{gd}} h_x, \quad \text{on } y=0\\
    \psi_x(x,y,0) = 0&,\quad \psi_y(x,y,0)=0,
    \end{rcases}
\end{align}
which only involves $\psi$ and $h$. We can then finally recover $f$ by the following procedure: We know by the fourth  and fifth equation in \eqref{reduced eq psi} that 
\begin{align*}
    f_x = p_x = -\psi_{y t} + C \psi_{x y} - \sqrt{\frac{d}{g}}A \psi_x
\end{align*}
on $y = 1$ and hence
\begin{align}\label{t derivative of f}
    f_t = - \psi_x- C f_x = C\psi_{y t} - C^2 \psi_{x y} + \left(\sqrt{\frac{d}{g}}AC -1\right)\psi_x
\end{align}
on $y = 1$.
Together with the initial condition
\begin{align*}
    f(x,0) = 0
\end{align*}
we can recover $f$ by
\begin{align*}
    f(x,t) = \int_{0}^t f_t(x,\tau)d \tau\,,
\end{align*}
and $p$ can be reconstructed similarly.

\section{General solution formulae for linear waves}

We will consider the space and space-time Fourier transform with notations
\begin{align*}
    \hat{\varphi}(\xi, y , t ) = \frac{1}{\sqrt{2 \pi}}\int_\R \varphi(x,y,t)e^{- i x \xi } dx
\end{align*}
and 
\begin{align*}
    \tilde{\varphi}(\xi,y,\omega) = \frac{1}{2 \pi} \int_\R \int_\R \varphi(x,y,t) e^{- i (x \xi + t\omega)}dx dt.
\end{align*}
We also define Fourier multipliers: Let $m :\R \rightarrow \cc$ be some function. We define $m(D)$ by 
\begin{align*}
    \left(m(D)\varphi\right)(x) = \frac{1}{\sqrt{2 \pi}}\int_\R m(\xi)\hat{\varphi}(\xi)e^{i x \xi }d\xi,
\end{align*}
or equivalently
\begin{align*}
    \left(\widehat{m(D)\varphi}\right)(\xi) = m(\xi) \hat{\varphi}(\xi),
\end{align*}
where $D = - i\partial_x$, whenever this is defined. Note that $m(D)$ maps real-valued functions to real-valued functions whenever the condition 
\begin{align}\label{FM condition}
    m(-\xi) = \overline{m(\xi)}
\end{align}
is satisfied. In particular, this is satisfied if $m$ is real-valued and even.
We can apply the space-time Fourier transform to the first three equations of system \eqref{reduced eq final} to obtain
\begin{align}\label{FT of main eq}
    \begin{rcases}
    \tilde{\psi}_{yy} - \delta^2\xi^2 \tilde{\psi} &= 0,\quad \text{in } 0< y < 1\\
    \left(1 + \sqrt{\frac{d}{g}}A C\right)\xi^2\tilde{\psi} &= Q^2 \tilde{\psi}_y - \sqrt{\frac{d}{g}}A \xi\omega\tilde{\psi},\quad \text{on }y=1\\
    \xi\tilde{\psi} &= - \omega  \tilde{h} - \frac{B}{\sqrt{gd}} \xi\tilde{h}, \quad \text{on } y=0
    \end{rcases},
\end{align}
where 
\begin{align}
    Q(\xi,\omega) = C \xi + \omega
\end{align}
with the relation
\begin{align}\label{T and Q}
    \widetilde{S\phi} = i Q\Tilde{\phi}.
\end{align}
Since \eqref{FT of main eq} is just a linear ODE of second order in $y$ for fixed $\xi$ and $\omega$, we can, somewhat explicitly, write down the solution:
\begin{align}\label{FT of psi by D1 and D2}
    \tilde{\psi}(\xi,y,\omega) = D_1(\xi,\omega)e^{\delta \xi y} + D_2(\xi,\omega)e^{-\delta \xi y},
\end{align}
where $D_1,D_2$ are determined by the boundary conditions. We can then reconstruct $\tilde{f}$ by applying the Fourier transform on the first equation in \eqref{t derivative of f} to obtain
\begin{align*}
    \omega \tilde{f}(\xi,\omega) = - \xi \tilde{\psi}(\xi,1,\omega) - C \xi \tilde{f}(\xi,\omega)
\end{align*}
and hence
\begin{align}\label{FT of f by FT of psi}
    \tilde{f}(\xi,\omega) = \frac{-\xi}{ C\xi+\omega}\tilde{\psi}(\xi,1,\omega) = \frac{-\xi}{Q(\xi,\omega)}\tilde{\psi}(\xi,1,\omega).
\end{align}
We ignore the potential singularity for the time being.
We know that the functions $D_1$ and $D_2$ are uniquely determined. Once we have formulae for them, we get a formula for $\tilde{\psi}$ by inserting in \eqref{FT of psi by D1 and D2}. Then we can insert this formula in \eqref{FT of f by FT of psi} to obtain an expression for $\tilde{f}$ depending on $\tilde{h}$.

We insert \eqref{FT of psi by D1 and D2} in \eqref{FT of main eq} to obtain 
\begin{align*}
    \xi\left(D_1(\xi,\omega) + D_2(\xi,\omega)\right) = -\left(\omega + \frac{B}{\sqrt{g d}}\xi\right)\Tilde{h}(\xi,\omega)
\end{align*}
and
\begin{align*}
    \left(\xi + \sqrt{\frac{d}{g}}A Q(\xi,\omega)\right)\left(D_1(\xi,\omega)e^{\delta \xi} + D_2(\xi,\omega)e^{-\delta \xi}\right) \\= \delta  Q(\xi,\omega)^2 \left(D_1(\xi,\omega) e^{\delta \xi} - D_2(\xi,\omega) e^{-\delta \xi}\right).
\end{align*}
This linear equation in $D_1(\xi,\omega)$ and $D_2(\xi,\omega)$ is solved easily and by \eqref{FT of f by FT of psi} we have 
\begin{align*}
    \Tilde{f}(\xi,\omega) &= \frac{-\xi}{Q(\xi,\omega)}\Tilde{\psi}(\xi,1,\omega))\\
    &=\frac{\delta  Q(\xi,\omega)\left(\omega + \frac{B}{\sqrt{g d}}\xi\right)}{\left(\delta Q(\xi,\omega)^2\cosh(\delta \xi) - \left(\xi + \sqrt{\frac{d}{g}}A Q(\xi,\omega)\right)\sinh(\delta \xi)\right)}\Tilde{h}(\xi,\omega),
\end{align*}
which we will write as
\begin{align}\label{FT of f by FT of h}
    \Tilde{f}(\xi,\omega) = \frac{\delta^2  Q(\xi,\omega)\left(\omega + \frac{B}{\sqrt{g d}}\xi\right)\Tilde{h}(\xi,\omega)}{\left(\delta^2 Q(\xi,\omega)^2\cosh(\delta \xi) - \left(\delta\xi + \sqrt{\frac{d}{g}}A \delta Q(\xi,\omega)\right)\sinh(\delta \xi)\right)}.
\end{align}
We will also write this as
\begin{align}\label{f by h with multipliers}
    \left( S^2 + \frac{D -i \sqrt{\frac{d}{g}}A S}{\delta}\tanh(\delta D)\right)f = \frac{  S\left(\partial_t + i\frac{B}{\sqrt{g d}}D\right)}{\cosh(\delta D)}h.
\end{align}
One would need to consider the roots of the denominator in \eqref{FT of f by FT of h} in order to justify \eqref{f by h with multipliers} rigorously. However for the case $A=0$, it turns out that the singularity cancels.

In the case $A = 0$ formula \eqref{f by h with multipliers} has the simpler form 
\begin{align}\label{f by h with multipliers, A=0}
    \left(S^2 + \frac{D}{\delta}\tanh(\delta D)\right)f = \frac{S^2}{\cosh(\delta D)}h.
\end{align}
We will extract a formula for $f$ whenever it solves the following type of equation. 
Let $\tau:\R \rightarrow \R$ be an even function such that $\tau(\xi) > 0$ for all $\xi\not=0$. Note that in particular \eqref{FM condition} is satisfied. We try to derive a formula for the solution $f$ of
\begin{align}\label{f by theta with multipliers}
    \begin{rcases}
    \left(S^2 + \tau(D)\right)f &= \theta\\
    \lim_{t \to -\infty} f &= 0
    \end{rcases},
\end{align}
where $\theta$ is a given forcing term. We will assume $f$ and $\theta$ to be in $\S(\R^2,\R)$ and the limit in \eqref{f by theta with multipliers} is adequate for functions in the Schwartz class \cite{Str94}. We can decompose the operator $S^2 + \tau(D)$ into
\begin{align}\label{factorizing T and tau}
    S^2 + \tau(D) = (S + i\sqrt{\tau(D)})(S - i\sqrt{\tau(D)}),
\end{align}
since $S$ and $i\sqrt{\tau(D)}$ commute. If we now set
\begin{align*}
    f_1 = (S + i\sqrt{\tau(D)})f
\end{align*}
then we see by \eqref{factorizing T and tau} and the fact that we can commute the two factors that
\begin{align*}
    (S - i\sqrt{\tau(D)})f_1 = \theta.
\end{align*}
We can now filter $f_1$ by the group $e^{i t \sqrt{\tau(D)}}$: Let $f_2 = e^{- i t \sqrt{\tau(D)}}f_1$. A calculation shows that 
\begin{align*}
    Sf_2 = (C\partial_x + \partial_t)f_2 = e^{- i t \sqrt{\tau(D)}}\theta.
\end{align*}
This is an inhomogeneous transport equation with solution
\begin{align*}
    f_2(x,t) = \int_{- \infty}^t e^{- i s \sqrt{\tau(D)}}\theta(x + C(s - t),s)ds.
\end{align*}
This yields
\begin{align*}
    \left(\sqrt{\tau(D)}f\right)(x,t) = \Imag \left(\int_{- \infty}^t e^{ i (t-s) \sqrt{\tau(D)}}\theta(x + C(s - t),s)ds\right).
\end{align*}
If we write out the Fourier multipliers explicitly, we get 
\begin{align}\label{explicit FM, A=0}
    f(x,t) = \frac{1}{\sqrt{2 \pi}} \Imag\left(\int_{-\infty}^t \int_\R  \frac{e^{i(t-s)\sqrt{\tau(\xi)}}}{\sqrt{\tau(\xi)}}e^{i x \xi} e^{iC(s-t)\xi} \hat{\theta}(\xi,s)d\xi ds\right).
\end{align}
We will deal with the possible singularity in the for us relevant case:
We have $\tau(\xi) = \frac{\xi}{\delta}\tanh(\delta \xi)$ and $\theta = \frac{S^2}{\cosh(\delta D)}h$. We see that $\tau$ is even and $\tau(\xi) >0$ for $\xi \not=0$. We have furthermore that $\tau$ is smooth with derivatives of at most polynomial growth at infinity. We will use that 
\begin{align}\label{limit of tau near zero}
    \lim_{\xi \to 0}\frac{|\xi|}{\sqrt{\tau(\xi)}}= 1.
\end{align}
Now observe that the function $\frac{1}{\sqrt{\tau(D)}}\theta$ is real-valued and hence we can subtract \\ $\left(\frac{1}{\sqrt{\tau(D)}}\theta\right)(0,\cdot)$ in \eqref{explicit FM, A=0} to obtain
\begin{align}\label{formula for f with cancelled sing, A=0}
    f(x,t) = \frac{1}{\sqrt{2 \pi}} \Imag\left(\int_{-\infty}^t \int_\R  \frac{e^{i(t-s)\sqrt{\tau(\xi)}}e^{i x \xi} e^{iC(s-t)\xi} - 1}{\sqrt{\tau(\xi)}} \hat{\theta}(\xi,s)d\xi ds\right).
\end{align}
Clearly, $\hat{\theta}(\xi,s)$ remains bounded around 0. For the quotient, note that we can rewrite it as follows:
\begin{align*}
    \left|\frac{e^{i(t-s)\sqrt{\tau(\xi)}}e^{i x \xi} e^{iC(s-t)\xi} - 1}{\sqrt{\tau(\xi)}} \right|
    =\left|\frac{e^{i(t-s)\sqrt{\tau(\xi)}}e^{i x \xi} e^{iC(s-t)\xi} - 1}{\xi}\right|\left|\frac{\xi}{\sqrt{\tau(\xi)}}\right| 
\end{align*}
The second factor converges to 1 by \eqref{limit of tau near zero}. For the first factor, note that the function
\begin{align*}
    \xi \mapsto e^{i(t-s)\sqrt{\tau(\xi)}}e^{i x \xi} e^{iC(s-t)\xi}
\end{align*}
is left and right differentiable in 0 and hence the first factor remains bounded as $\xi \to 0$ for any $s$. We conclude that \eqref{formula for f with cancelled sing, A=0} holds (at least pointwise).

Note that we cannot do the same for the case $A\not=0$: The operator on the left hand side of \eqref{f by h with multipliers} is given by 
\begin{align*}
    S^2 + (1 - C^2)\frac{D}{\delta}\tanh(\delta D) - \frac{iC}{\delta}\partial_t \tanh(\delta D)
\end{align*}
and there does not seem to be an apparent way to obtain a root for
\begin{align*}
    (1 - C^2)\frac{D}{\delta}\tanh(\delta D) -\frac{iC}{\delta}\partial_t \tanh(\delta D).
\end{align*}
In fact, it is unclear whether this operator is even positive.

\section{Behaviour in different regimes}

We want to investigate some regimes for which our considerations give useful predictions. We will consider on the one hand the case where $\delta$ is small, which is justified considering the model as $\delta \to 0$. On the other hand, we will discuss the case where $\delta$ is finite and non-vanishing. The case where $\delta$ is large (corresponding to the formal limit $\delta \to \infty$) is not relevant for us since deep-water tsunamis are rare.
\subsection{The shallow water regime}
For simplicity, we use the model in the case $A=0$: As $\delta \to 0$, \eqref{f by h with multipliers, A=0} becomes 
\begin{align}\label{shallow water pde}
    (S^2 - \partial_x^2)f = S^2 h.
\end{align}
The initial conditions are on the one hand $f(x,0)=0$ from \eqref{reduced eq}. On the other hand, this gives $f_x(x,0)=0$ and inserting this in the sixth relation in \eqref{reduced eq} for $t=0$ yields $f_t(x,0) = 0$. By using Duhamel's principle (similar to \cite[p. 80--81]{Eva10}), we get the formula
\begin{align}\label{formula for f in shallow water}
    f(x,t) = \frac{1}{2} \int_0^t \int_{x - (t-s)(C+1)}^{x - (t-s)(C-1)}S^2 h(r,s) dr ds.
\end{align}

We will assume that $h$ can be separated in the following way: $h(x,t) = a(t)b(x)$ for $a \in C^2(\R,[0,\infty))$, $b \in C^2(\R,\R)$ such that $a(t)=0$ for $t \leq 0$ and $a(t) = 1$ for $t > t_0$ (where $t_0$ represents the duration of the earthquake), and with $b(x) = 0$ for $x \notin (-L,L)$ modelling a localized tsunami source. Inserting this into \eqref{formula for f in shallow water} yields the formula
\begin{align*}
    f(x,t)  &=   a(t)b(x) \\
            &+\frac{1}{2}\int_0^t a(s)\left[b'(x - (t-s)(C-1)) - b'(x - (t-s)(C+1))\right]ds,
\end{align*}
which is a lot simpler. In the limiting case $t_0 \searrow 0$, $a$ becomes the Heaviside step function (modelling an instantaneous upward thrust of the seabed near to the earthquake's epicentre) and the formula further simplifies to
\begin{align}\label{final formula for shallow water}
    f(x,t)  =\frac{C^2}{C^2-1}b(x) + \frac{b(x - t(1+C))}{2(1+C)} + \frac{b(x + t(1-C))}{2(1-C)},
\end{align}
for $x \in \R$ and $t >0$. We can draw some insightful conclusions:
\begin{itemize}
    \item at each instant $t >0$ after initiation, the generated wave is localized as $f(x,t) = 0$ for $|x| \geq L + t(1+C)$;
    \item the surface wave consists of one stationary part (which can be disregarded, since $C \ll 1$) and two travelling waves, one moving to the right and the other moving to the left;
    \item the wave travelling to the left moves with non-dimensionalized speed $1-C$ (corresponding to $\sqrt{gd} - B$ in the physical variables) and the wave travelling to the right moves with non-dimensionalized speed $1+C$ (corresponding to $\sqrt{gd} + B$ in the physical variables);
    \item the shapes of the waves travelling to the left and to the right remain unchanged and are precisely that of the bed deformation at a scale of $\frac{1}{2(1-C)}$ and $\frac{1}{2(1+C)}$, respectively.
\end{itemize}

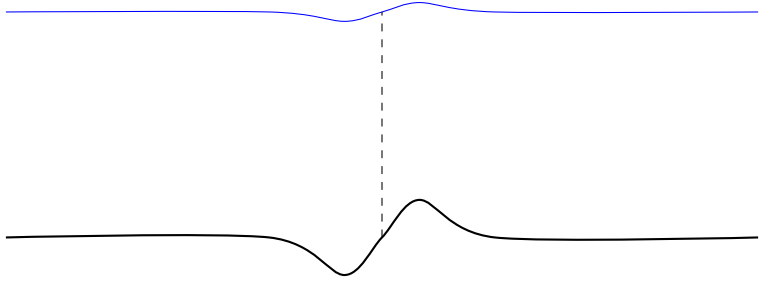
\begin{figure}[!ht]
\centering
\begin{tikzpicture}
\draw[blue] plot[smooth] coordinates {(-5,3) (-1.5,3) (-0.5,2.875) (0,3) (0.5,3.125) (1.5,3) (5,3)};
\draw[thick] plot[smooth] coordinates {(-5,0) (-1.5,0) (-0.5,-0.5) (0,0) (0.5,0.5) (1.5,0) (5,0)}  ;
\draw[dashed] (0,0) -- (0,3);
\end{tikzpicture}
\caption{Model shortly after the instantaneous upward and downward thrust}
\label{Model at time t small}
\end{figure}

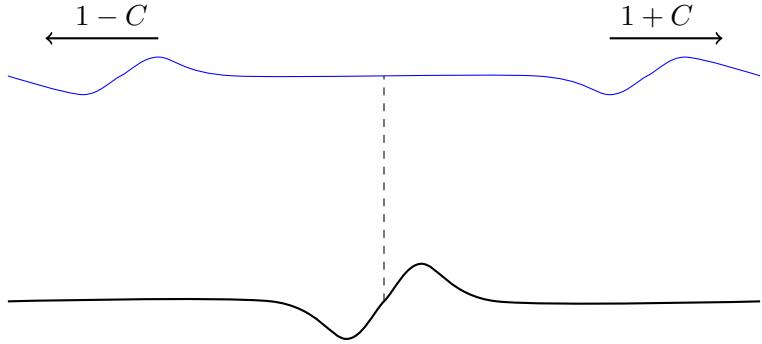
\begin{figure}[!ht]
\centering
\begin{tikzpicture}
\draw[blue] plot[smooth] coordinates {(-5,3) (-4,2.75) (-3.5,3) (-3,3.25) (-2,3) (2,3) (3,2.75) (3.5,3) (4,3.25) (5,3)};
\draw[thick] plot[smooth] coordinates {(-5,0) (-1.5,0) (-0.5,-0.5) (0,0) (0.5,0.5) (1.5,0) (5,0)}  ;
\draw[dashed] (0,0) -- (0,3);
\draw[thick, <-] (4.5,3.5) -- (3,3.5) node[anchor = south west] {$1 + C$};
\draw[thick, <-] (-4.5,3.5) -- (-3,3.5) node[anchor = south east] {$1 - C$};
\end{tikzpicture}
\caption{Expected long-term behaviour}
\label{Model at time t large}
\end{figure}
This behaviour is illustrated in Figure \ref{Model at time t small} and Figure \ref{Model at time t large}.
On the one hand, the wave travelling to the right is faster for bigger $C$, but the maximal amplitude decreases. On the other hand, the wave travelling to the left becomes slower for bigger $C$, but the maximal amplitude increases.

For $B=0$ and hence $C=0$, this model gives good predictions for the maximal amplitude and speed of the propagating wave for the two largest historical tsunamis: the December 2004 and the May 1960 tsunamis, see \cite{CG12}. Since $B \ll \sqrt{gd}$, i.e. $C \ll 1$, the model for $B \not= 0$ will give a very similar prediction, which might be slightly more precise, if one were in a situation, where $B$ is known somewhat precisely, say due to a current.

\subsection{The finite, non-vanishing regime}

Here we cannot expect to be able to write down $f$ in any explicit fashion. However, we will discuss a method to analyze the asymptotic behaviour of $f$ in \eqref{formula for f with cancelled sing, A=0}. Ultimately, we will see that this approach does not work, since the method only works for too large $t$, see \eqref{stationary-phase analysis takes too long}.

We cannot do the same for $B = 0$ or the general case, since we lack a formula which would need to correspond to \eqref{explicit FM, A=0} or \eqref{formula for f with cancelled sing, A=0}. If such a formula were available, stationary-phase analysis might then be appropriate.
\subsubsection{Principle of stationary-phase}
We want to analyse the asymptotic behaviour of \eqref{formula for f with cancelled sing, A=0}. For this we try to use the stationary-phase principle.
Recall that we set $\tau(\xi) = \frac{\xi}{\delta}\tanh(\delta \xi)$ and $\theta = \frac{S^2}{\cosh(\delta D)}h$.
We first rewrite \eqref{formula for f with cancelled sing, A=0} as
\begin{align*}
    &f(x,t) = \frac{1}{\sqrt{2 \pi}} \Imag\left(\int_{-\infty}^t \int_\R  \frac{e^{i(t-s)\sqrt{\tau(\xi)}}e^{i x \xi} e^{iC(s-t)\xi} - 1}{\sqrt{\tau(\xi)}} \hat{\theta}(\xi,s)d\xi ds\right) \\
            &= \frac{1}{\sqrt{2 \pi}} \Imag\left(\int_{-\infty}^t \int_\R  \frac{e^{-i s \left( \sqrt{\tau(\xi)} - C \xi\right)}e^{it\left(\sqrt{\tau(\xi)} + (\mathcal{X}-C)\xi)\right)} - 1}{\sqrt{\tau(\xi)}} \hat{\theta}(\xi,s)   d\xi ds\right),
\end{align*}
where $\X = \frac{x}{t}$. The derivative of the phase factor $\sqrt{\tau(\xi)} + (\mathcal{X}-C)\xi$ is given by 
\begin{align*}
    [\sqrt{\tau}]'(\xi) + \X - C
\end{align*}
and one can check that 
\begin{align}\label{first derivative of root tau}
    \partial_\xi\sqrt{\tau(\xi)}=\frac{\sinh(\delta \xi)\cosh(\delta \xi) + \delta \xi}{2\delta\cosh^2(\delta \xi)} \sqrt{\frac{\delta\cosh(\delta \xi)}{\xi\sinh(\delta \xi)}}
\end{align}
and
\begin{align}\label{second derivative of root tau}
    \partial_\xi^2\sqrt{\tau(\xi)} = -\frac{1}{\tau^\frac{3}{2}(\xi)}\frac{\left(\sinh(\delta \xi)\cosh(\delta \xi) - \delta \xi\right)^2 + 4\delta^2\xi^2 \sinh^2(\delta \xi)}{4 \delta^2 \cosh^4(\delta \xi)} < 0.
\end{align}
One computes the limits
\begin{align*}
    &\lim_{\xi \searrow 0}\partial_\xi\sqrt{\tau(\xi)} = 1,\\
    &\lim_{\xi \nearrow 0}\partial_\xi\sqrt{\tau(\xi)} = -1,\\
    &\lim_{\xi \to \infty}\partial_\xi\sqrt{\tau(\xi)} = 0,\\
    &\lim_{\xi \to -\infty}\partial_\xi\sqrt{\tau(\xi)} = 0\,,
\end{align*}
and concludes that $\xi \mapsto \partial_\xi\sqrt{\tau(\xi)}$ is strictly decreasing on $(0,\infty)$ from the asymptotic value 1 towards the asymptotic value 0 and on $(-\infty,0)$ from the asymptotic value $0$ to the asymptotic value $-1$. So we conclude that 
\begin{align*}
    [\sqrt{\tau}]'(\xi_0) + \X - C = 0
\end{align*}
is possible if $\X - C \in (-1,0)\cup(0,1)$ for exactly one $\xi_0$ for given $\X$. So as long as $\X$ is in the corresponding range, we can write $\xi_0 = \xi_0(\X)$.
The stationary-phase principle gives
\begin{align*}
    f(x,t)\sim &\frac{1}{\sqrt{t \tau(\xi_0) |[\sqrt{\tau}]''(\xi_0)|}}\cdot \\&\Imag\left(\int_{-\infty}^t  e^{-i s \left( \sqrt{\tau(\xi_0)} - C \xi_0\right)} \hat{\theta}(\xi_0,s) e^{it\left(\sqrt{\tau(\xi_0)} + (\mathcal{X}-C)\xi_0 + \sigma\frac{\pi}{4})\right)}   ds\right)
\end{align*}
for large $t$, where $\sigma$ is the sign of $[\sqrt{\tau}]''(\xi_0)$, which is just $-1$ by \eqref{second derivative of root tau}. We see that this is a Fourier transform (with respect to time) and hence get
\begin{align*}
    f(x,t) \sim &\frac{\sqrt{2\pi}}{\sqrt{t \tau(\xi_0) |[\sqrt{\tau}]''(\xi_0)|}}\cdot \\&\Imag\left(  \Tilde{\theta}(\xi_0,\sqrt{\tau(\xi_0)} - C \xi_0) e^{it\left(\sqrt{\tau(\xi_0)} + (\mathcal{X}-C)\xi_0 -\frac{\pi}{4})\right)}   \right).
\end{align*}
We have
\begin{align*}
    \Tilde{\theta}(\xi_0,\sqrt{\tau(\xi_0)} - C \xi_0) &= -\frac{Q(\xi_0,\sqrt{\tau(\xi_0)} - C \xi_0)^2}{\cosh(\delta \xi_0)}\Tilde{h}(\xi_0,\sqrt{\tau(\xi_0)} - C \xi_0) \\
    &=-\frac{\tau(\xi_0)}{\cosh(\delta \xi_0)}\Tilde{h}(\xi_0,\sqrt{\tau(\xi_0)} - C \xi_0)
\end{align*}
and hence get 
\begin{align}\label{f stationary-phase asymptotic}
    f(x,t) \sim &\frac{-\sqrt{2\pi \tau(\xi_0)}}{\cosh(\delta \xi_0)\sqrt{t  |[\sqrt{\tau}]''(\xi_0)|}}\cdot \nonumber\\&\Imag\left(   e^{it\left(\sqrt{\tau(\xi_0)} + (\mathcal{X}-C)\xi_0 -\frac{\pi}{4})\right)} \Tilde{h}(\xi_0,\sqrt{\tau(\xi_0)} - C \xi_0)  \right).
\end{align}
Finally, we try to extract the asymptotics with respect to $\delta$: We fix a $\xi > 0$. A look at \eqref{first derivative of root tau} gives that $\delta \mapsto \partial_\xi \sqrt{\tau(\xi)}$ is $C^\infty$ and clearly even in $\delta$. Since we have $\lim_{\delta \searrow 0}\partial_\xi \sqrt{\tau(\xi)} = 1$, the Taylor expansion with respect to $\delta$ around 0 has the form
\begin{align*}
    \partial_\xi \sqrt{\tau(\xi)} = 1 + \alpha(\xi)\delta^2 + O(\delta^3),
\end{align*}
for some coefficient $\alpha(\xi)$. One could find this coefficient through tedious calculations, but one can instead notice that $2 \alpha(\xi)$ is the coefficient of $\delta^2$ of the expansion of $\left(\partial_\xi \sqrt{\tau(\xi)}\right)^2$. We have
\begin{align*}
    \left(\partial_\xi \sqrt{\tau(\xi)}\right)^2 &=\frac{\sinh(\delta \xi)}{4\delta \xi\cosh(\delta \xi)} + \frac{1}{2\cosh^2(\delta \xi)} + \frac{\delta \xi}{4\sinh(\delta \xi)\cosh^3(\delta \xi)}\\
    &= 1 - \delta^2 \xi^2 + O(\delta^3)
\end{align*}
and hence we conclude $\partial_\xi \sqrt{\tau(\xi)} = 1 - \frac{1}{2}\delta^2\xi^2 + O(\delta^3)$. We conclude that for $\X - C \in (-1,0)$ the stationary-phase point $\xi_0 >0$ with $[\sqrt{\tau}]'(\xi_0) + \X - C = 0$ satisfies $\delta \xi_0 = O(1)$. Inserting this into \eqref{second derivative of root tau} yields 
\begin{align}\label{second derivative at station}
    |\partial_\xi^2\sqrt{\tau(\xi_0)}| = O(\delta).
\end{align}
\subsubsection{Some typical physical parameters}

We want to see whether stationary-phase analysis, using the formula 
\begin{align*}
    f(x,t) \sim &\frac{-\sqrt{2\pi \tau(\xi_0)}}{\cosh(\delta \xi_0)\sqrt{t  |[\sqrt{\tau}]''(\xi_0)|}}\cdot \\&\Imag\left(   e^{it\left(\sqrt{\tau(\xi_0)} + (\mathcal{X}-C)\xi_0 -\frac{\pi}{4})\right)} \Tilde{h}(\xi_0,\sqrt{\tau(\xi_0)} - C \xi_0)  \right),
\end{align*}
could be justified here. Since the stationary-phase principle applies for large $t$ we insert typical values of the physical parameters for tsunamis propagating at open sea (see \cite{Con09}): 
\begin{align*}
    a = 1\text{m},\quad d=4 \text{km},\quad \lambda = 200 \text{km}.
\end{align*}
This leads to $\eps = 0.00025$ and $\delta = 0.02$ which point to linear wave approximation. Applying the stationary-phase principle is justified if 
\begin{align}
    t |\sqrt{\tau}''(\xi)| \gg 1.
\end{align}
We have by \eqref{second derivative at station}, $|\sqrt{\tau}''(\xi)| = O(\delta)$, which translates to 
\begin{align}\label{t gg delta squared}
    t \gg \frac{1}{\delta^2}.
\end{align}
In physical variables this means
\begin{align}\label{stationary-phase analysis takes too long}
    T \gg \frac{\lambda}{\sqrt{gd}\delta^2} \approx 2.5 \times 10^6 \text{ s} \approx 700 \text{ h},
\end{align}
which takes way too long to be applicable.
If one were to relax \eqref{t gg delta squared} to a condition of the form
\begin{align*}
    t \geq \frac{1}{\delta},
\end{align*}
one would arrive at 
\begin{align*}
    T \geq \frac{\lambda}{\sqrt{gd}\delta} \approx 5 \times 10^4 \text{ s} \approx 14 \text{ h}
\end{align*}
which would barely be applicable in the case where the epicentre of the seaquake is far off the shore. Nevertheless, relaxing \eqref{t gg delta squared} is quite dubious. 
We conclude that using stationary-phase analysis is not an effective tool in the case $A = 0$.

\end{document}